\documentclass[a4paper,11pt]{amsart}
\usepackage{fancyhdr}
\usepackage{amsmath}
\usepackage{dsfont}
\usepackage{hyperref}
\usepackage[mathscr]{eucal}
\usepackage[cp1251]{inputenc}
\usepackage[english]{babel}
\usepackage{enumerate,float,indentfirst}
\usepackage{graphicx}
\usepackage{xcolor}
\usepackage{latexsym,a4,mathrsfs,amsthm,amsmath,amssymb,url}
\usepackage{amsfonts}
\usepackage{amssymb}

\usepackage{tikz}
\usepackage{pgfplots}
\usetikzlibrary{plotmarks}

\numberwithin{equation}{section}
\newcommand{\mns}{\scalebox{0.35}[1.0]{\( - \)}}
\newcommand{\dist}{\mathop{\rm dist}\nolimits}
\newcommand{\Span}{\mathop{\rm Span}\nolimits}

\begin{document}

\vspace{1in}

\title[Spectral boundary conditions for volumetric frame fields design]{\bf Spectral boundary conditions \\ for volumetric frame fields design}

\author[Yu. Nesterenko]{Yu. Nesterenko}
\address{ Siemens Digital Industries Software }
\email{Yuri.Nesterenko@siemens.com}

\begin{abstract}

Using the 4th and the 3rd degree spherical harmonics as the representations for volumetric frames, we describe a simple algebraic technique for combining multiple frame orientation constraints into a single quadratic penalty function. This technique allows to solve volumetric frame fields design problems using a coarse-to-fine strategy on hierarchical grids with immersed boundaries. These results were presented for the first time at the FRAMES 2023 European workshop on meshing.

\end{abstract}

\maketitle

\thispagestyle{empty}

\section{Introduction}

In this research note we consider the question: how to enforce multiple boundary conditions in a volumetric frame fields design problem? The reasons of such multiplicity may be different: noisy input data, presence of multiple orientation constraints (e.g. in a CAD model with a piecewise smooth surface), or consideration of boundary conditions in an immersed way. The last scenario, in our opinion, is of particular interest, since it opens up the possibility to optimize volumetric frame fields using a coarse-to-fine strategy on hierarchical grids with immersed boundaries.

Based on the two different representations of volumetric frames --- the 4th and the 3rd degree spherical harmonics --- we propose a simple algebraic technique allowing to combine arbitrary many frame orientation constraints into a single quadratic penalty function (of 9 and 7 variables respectively).

\section{Spherical harmonics of degree 4}

We consider real-valued spherical harmonics of 4th degree on the unit sphere.
These functions form 9D vector space with the standard orthonormal basis $Y_{4,-4}, \ldots, Y_{4,4}$ (see \cite{Gorller1996}).

\begin{figure}[h]
\includegraphics[width=13.0cm]{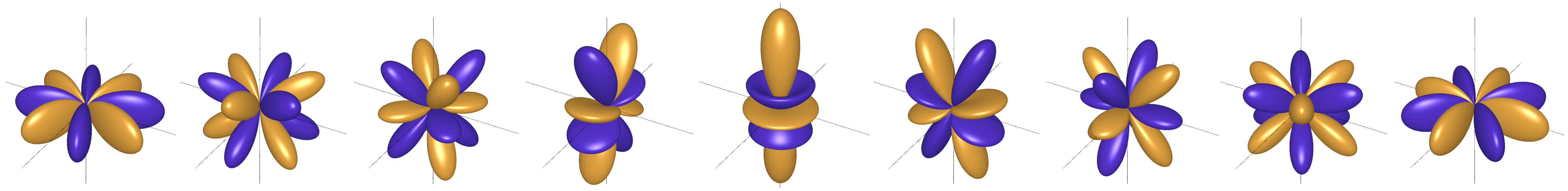}
\caption{ Spherical plots of basis functions $Y_{4,\mns4}, \ldots, Y_{4,4}$. }
\label{fig:Y4i}
\end{figure}

In the coordinate form, all octahedrally symmetric spherical harmonics (up to multiplication by $-1$) may be obtained from the reference one ---
\begin{equation*}
\tilde{h} = (0,0,0,0,\sqrt{\frac{7}{12}},0,0,0,\sqrt{\frac{5}{12}})^T \in \mathds{R}^9
\end{equation*}
--- by rotations
\begin{equation*}
h = R_x(\alpha) \times R_y(\beta) \times R_z(\gamma) \times \tilde{h},
\end{equation*}
where $\alpha, \beta, \gamma$ are Euler angles.

\begin{figure}[h]
\includegraphics[width=8.0cm]{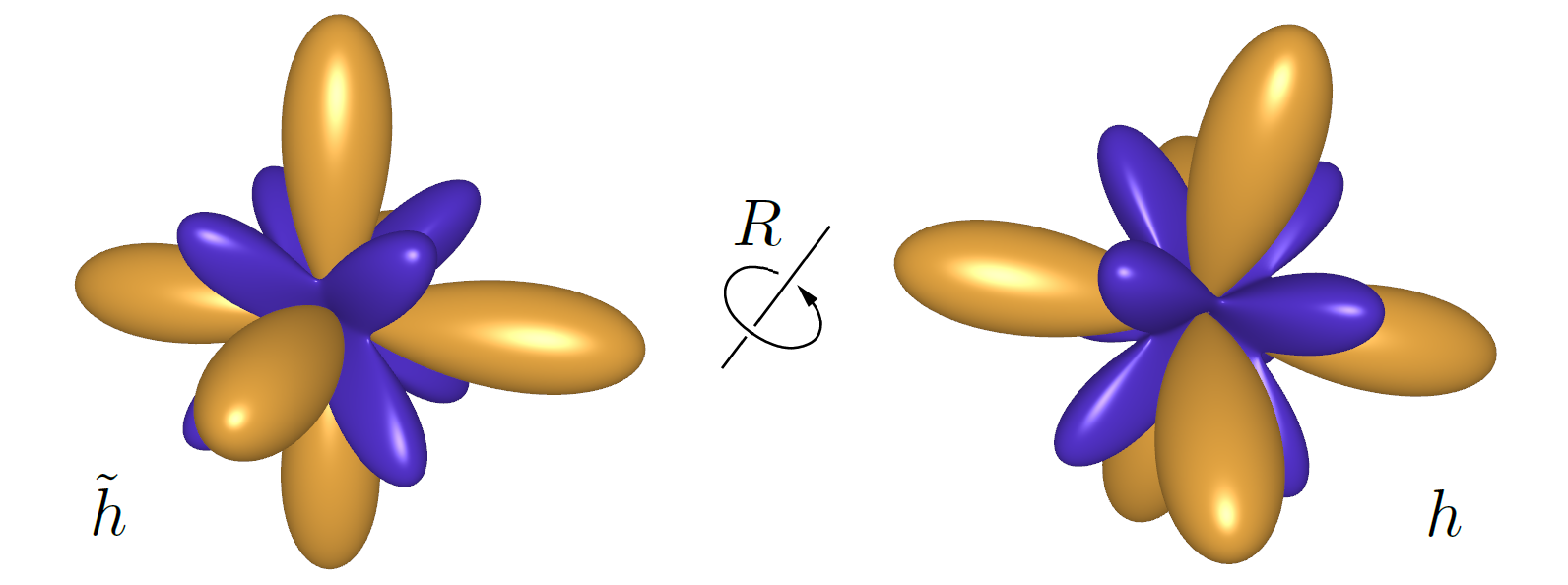}
\caption{ The reference harmonic and its rotation. }
\label{fig:def}
\end{figure}
Appendix A.1 describes the construction of the rotation matrices $R_x$, $R_y$ and $R_z$.

If we restrict ourselves to rotations about $z$ axis, we get the next 1D manifold of spherical harmonics:
\begin{equation*}
\begin{split}
& h = (0,0,0,0,\sqrt{\frac{7}{12}},0,0,0,0)^T + \\
\alpha \, (1,0,0,& 0,0,0,0,0,0)^T + \beta \, (0,0,0,0,0,0,0,0,1)^T,
\end{split}
\end{equation*}
where $\alpha^2 + \beta^2 = 5/12$.

Further in this section, we will denote the coordinate vectors from the last formula as $b_z$, $p_z$ and $q_z$ respectively.

To enforce an arbitrary symmetric harmonic\footnote{We refer the reader to \cite{Nesterenko2020} for the definition of the symmetrization penalty term.} to lie on the given manifold --- and thus to respect the orientation constraints given by $z$ axis --- we can apply the penalty function
\begin{equation*}
E_z(h) = \dist^2 (h - b_z, \, \Span\{p_z, q_z\}).
\end{equation*} 

The similar formula holds for orientation constraints with an arbitrary normal vector $n$:
\begin{equation}\label{dist4}
E_n(h) = \dist^2 (h - b_n, \, \Span\{p_n, q_n\}),
\end{equation}
where $b_n = R_{z \rightarrow n} \times b_z$, $p_n = R_{z \rightarrow n} \times p_z$, $q_n = R_{z \rightarrow n} \times q_z$ and $R_{z \rightarrow n}$ is the corresponding $9 \times 9$ rotation matrix.

The formula (\ref{dist4}) may be rewritten as follows:
\begin{equation}\label{E4}
E_n(h) = h^T (I - p_n \, p_n^T - q_n \, q_n^T) \, h - 2 b_n^T \, h + b_n^T \, b_n,
\end{equation}
where the $9 \times 9$ matrix $I - p_n \, p_n^T - q_n \, q_n^T$ is the orthogonal projection onto $\Span^\perp\{p_n, q_n\}$ and the free term $b_n^T \, b_n$ is obviously equal to $\frac{7}{12}$.

Now we can enforce multiple orientation constraints by simply summing the corresponding coefficients of the quadratic penalty functions (\ref{E4}):
\begin{equation}\label{sE4}
\sum_n w_n \, E_n(h) = h^T \left(\sum_n w_n \, A_n\right) h - 2 \left(\sum_n w_n \, b_n\right)^T h + \sum_n w_n \, b_n^T \, b_n. 
\end{equation}
Here, $w_n$ are some positive weights (in case a weighted sum is needed) and $A_n = I - p_n \, p_n^T - q_n \, q_n^T$.

\section{Spherical harmonics of degree 3 (octupoles)}

In this section we will describe the technique similar to the discussed above but for the 7D space of octupoles --- real-valued spherical harmonics of degree 3.

Note that all notations from the previous section will be reintroduced here with the similar meaning.

\begin{figure}[h]
\includegraphics[width=13.0cm]{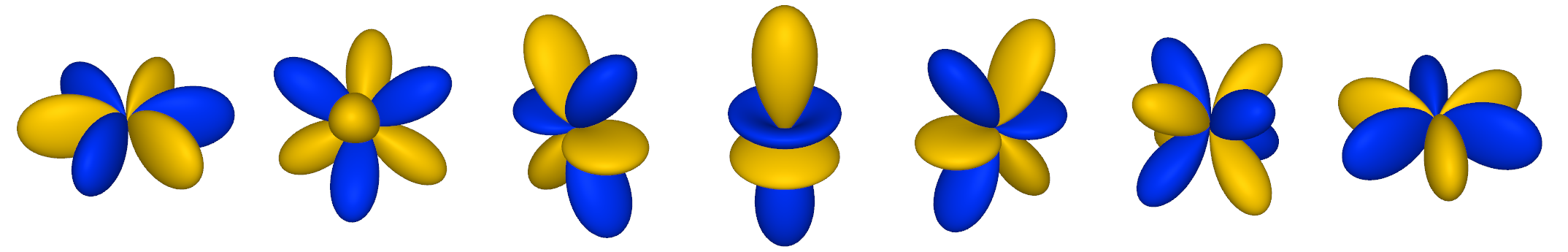}
\caption{ Spherical plots of basis functions $Y_{3,\mns3}, \ldots, Y_{3,3}$. }
\label{fig:Y3i}
\end{figure}

Let's consider the 3D manifold of \emph{semisymmetric} octupoles --- octupoles possessing octahedral symmetries up to multiplication by $-1$ (see \cite{Nesterenko2023}).

In the coordinate form, all octupoles of this kind may be obtained from the reference one ---
\begin{equation*}
\tilde{h} = (0,1,0,0,0,0,0)^T \in \mathds{R}^7
\end{equation*}
--- by rotations
\begin{equation*}
h = R_x(\alpha) \times R_y(\beta) \times R_z(\gamma) \times \tilde{h},
\end{equation*}
where $\alpha, \beta, \gamma$ are Euler angles.

\begin{figure}[h]
\includegraphics[width=8.0cm]{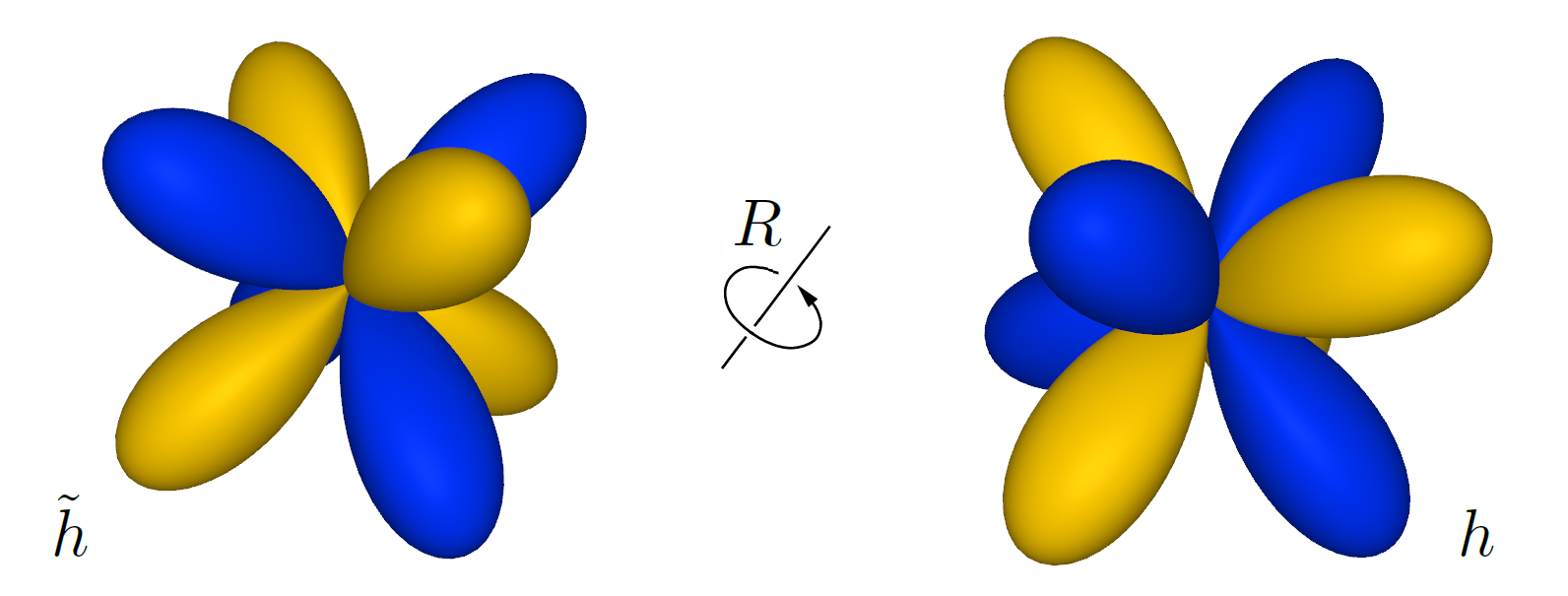}
\caption{ The reference octupole and its rotation. }
\label{fig:ref7}
\end{figure}
Appendix A.2 describes the construction of the rotation matrices $R_x$, $R_y$ and $R_z$ for the space of octupoles.

If we restrict ourselves to rotations about $z$ axis, we get the next 1D manifold of octupoles:
\begin{equation*}
h = \alpha \, (0,1,0,0,0,0,0)^T + \beta \, (0,0,0,0,0,1,0)^T,
\end{equation*}
where $\alpha^2 + \beta^2 = 1$.

Further, we will denote the coordinate vectors from the last formula as $p_z$ and $q_z$ respectively.

To enforce an arbitrary semisymmetric octupole\footnote{We refer the reader to \cite{Nesterenko2023} for the definition of the semisymmetrization penalty term.} to lie on the given 1D manifold --- and thus to respect the orientation constraints given by $z$ axis --- we can apply the penalty function
\begin{equation*}
E_z(h) = \dist^2 (h, \, \Span\{p_z, q_z\}).
\end{equation*} 

The similar formula holds for orientation constraints with an arbitrary normal vector $n$:
\begin{equation}\label{dist3}
E_n(h) = \dist^2 (h, \, \Span\{p_n, q_n\}),
\end{equation}
where $p_n = R_{z \rightarrow n} \times p_z$, $q_n = R_{z \rightarrow n} \times q_z$ and $R_{z \rightarrow n}$ is the corresponding $7 \times 7$ rotation matrix.

The formula (\ref{dist3}) may be rewritten as follows:
\begin{equation}\label{E3}
E_n(h) = h^T (I - p_n \, p_n^T - q_n \, q_n^T) \, h,
\end{equation}
where the $7 \times 7$ matrix $I - p_n \, p_n^T - q_n \, q_n^T$ is the orthogonal projection onto $\Span^\perp\{p_n, q_n\}$ subspace.

Now we can enforce multiple orientation constraints by simply summing the corresponding coefficients of the quadratic penalty functions (\ref{E3}):
\begin{equation}\label{sE3}
\sum_n w_n \, E_n(h) = h^T \left(\sum_n w_n \, A_n\right) h. 
\end{equation}
Here, $w_n$ are some positive weights (in case a weighted sum is needed) and $A_n = I - p_n \, p_n^T - q_n \, q_n^T$.

Note that in contrast to (\ref{sE4}), the quadratic penalty function (\ref{sE3}) is homogeneous and requires storing 28 coefficients of its symmetric $7\times7$ matrix instead of 55 coefficients of a symmetric $9\times9$ matrix together with a 9D vector and a constant (per boundary frame).

\section{Spectral analysis}

In this section we consider some particular combinations of orientation constraints and analyse their spectral properties. Among other things, the special role of spectra in the presented technique explains its name.

For simplicity, we will limit ourselves to
the octupoles case (using the corresponding notations).

\subsection{"Cube"}

Let's assume that the orientation constraints for some boundary octupole consist of the three mutually orthogonal axial-aligned normal vectors (taken with the unit weights).

The corresponding quadratic form has the matrix
\begin{equation*}
\begin{Small}
\begin{bmatrix}
\frac{21}{8}   & 0             & \frac{\sqrt{15}}{8} & 0             & 0             & 0             & 0             \\
0             & 0             & 0             & 0             & 0             & 0             & 0             \\
\frac{\sqrt{15}}{8} & 0             & \frac{19}{8}   & 0             & 0             & 0             & 0             \\
0             & 0             & 0             & 3             & 0             & 0             & 0             \\
0             & 0             & 0             & 0             & \frac{19}{8}   & 0             & \frac{-\sqrt{15}}{8} \\
0             & 0             & 0             & 0             & 0             & 2             & 0             \\
0             & 0             & 0             & 0             & \frac{-\sqrt{15}}{8} & 0             & \frac{21}{8}  
\end{bmatrix}
\end{Small}
\end{equation*}
with the eigenvalues $\{0, 2, 2, 2, 3, 3, 3\}$ and the null space $\Span\{ \tilde{h} \}$.

Thus, the given quadratic form penalizes all components of its argument, except the suitable one --- $\tilde{h}$.

\subsection{"Cylinder"}

If the boundary is a cylinder oriented along $y$ axis, then (up to a constant factor) the corresponding quadratic form has the matrix
\begin{equation*}
\frac{1}{2\pi} \int_{0}^{2\pi} (I - R_y(\beta) p_z \, p_z^T R_y(\beta)^T - R_y(\beta) q_z \, q_z^T R_y(\beta)^T) d\beta =
\end{equation*}
\begin{equation*}
\begin{Small}
\begin{bmatrix}
 \frac{13}{16} & 0 & \frac{\sqrt{15}}{16} & 0 & 0 & 0 & 0 \\
 0 & \frac{1}{2} & 0 & 0 & 0 & 0 & 0 \\
 \frac{\sqrt{15}}{16} & 0 & \frac{11}{16} & 0 & 0 & 0 & 0 \\
 0 & 0 & 0 & \frac{49}{64} & 0 & -\frac{\sqrt{15}}{64} & 0 \\
 0 & 0 & 0 & 0 & \frac{103}{128} & 0 & -\frac{\sqrt{15}}{128} \\
 0 & 0 & 0 & -\frac{\sqrt{15}}{64} & 0 & \frac{47}{64} & 0 \\
 0 & 0 & 0 & 0 & -\frac{\sqrt{15}}{128} & 0 & \frac{89}{128}
\end{bmatrix}
\end{Small}
\end{equation*}
with the eigenvalues $\{\frac{1}{2}, \frac{1}{2}, \frac{11}{16}, \frac{11}{16}, \frac{13}{16}, \frac{13}{16}, 1\}$.

It means that this quadratic form penalizes all the octupole components. To prevent such behavior we can simply shift the form spectra by subtracting $\frac{1}{2}I$ from its matrix. The null space of the shifted quadratic form --- $\Span\{\tilde{h}, R_y(\frac{\pi}{4}) \times \tilde{h}\}$ --- corresponds to the free rotation of the reference semisymmetric octupole $\tilde{h}$ about $y$ axis.

In other words, a thin cylinder immersed into a relatively coarse grid contributes to its cells' boundary conditions similarly to the axis of the given cylinder.

\subsection{"Shpere"}

Let's consider the unit sphere $S^2$.
It can be shown by direct calculations that
\begin{equation*}
\frac{1}{4 \pi} \int_{n \in S^2} A_n dS = \frac{5}{7} I.
\end{equation*}

Thus, taking into account the spectral shift, we can conclude that spheres lying entirely inside a grid cell don't contribute to its boundary conditions.  

\section{Conclusion}
Using the 4th and the 3rd degree spherical harmonics as the representations of volumetric frames, we have described the two types of quadratic penalty functions enforcing an arbitrary number of frame orientation constraints.

In combination with the (semi)symmetrization techniques from \cite{Nesterenko2020} and \cite{Nesterenko2023}, this approach allows to solve volumetric frame fields design problems using a coarse-to-fine strategy on hierarchical grids with immersed boundaries.

\section{Acknowledgements}

I would like to thank my friend Svyatoslav Proskurnya for a lot of fruitful discussions and, in particular, for giving the name to the presented technique long before it was formulated.

\section{Appendix A.1}
The rotational matrices $R_x$, $R_y$ and $R_z$ for spherical harmonics of degree 4 are defined as follows:

\begin{Small}
\begin{equation*}
R_z(\gamma) =
\begin{bmatrix}
\cos 4\gamma  & 0             & 0             & 0             & 0             & 0             & 0             & 0             &  \sin 4\gamma  \\
0             & \cos 3\gamma  & 0             & 0             & 0             & 0             & 0             &  \sin 3\gamma & 0              \\
0             & 0             & \cos 2\gamma  & 0             & 0             & 0             &  \sin 2\gamma & 0             & 0              \\
0             & 0             & 0             & \cos  \gamma  & 0             &  \sin  \gamma & 0             & 0             & 0              \\
0             & 0             & 0             & 0             & 1             & 0             & 0             & 0             & 0              \\
0             & 0             & 0             & -\sin  \gamma & 0             & \cos  \gamma  & 0             & 0             & 0              \\
0             & 0             & -\sin 2\gamma & 0             & 0             & 0             & \cos 2\gamma  & 0             & 0              \\
0             & -\sin 3\gamma & 0             & 0             & 0             & 0             & 0             & \cos 3\gamma  & 0              \\
-\sin 4\gamma & 0             & 0             & 0             & 0             & 0             & 0             & 0             & \cos 4\gamma
\end{bmatrix},
\end{equation*}

\begin{equation*}
R_x(\frac{\pi}{2}) = \frac{1}{8}
\begin{bmatrix}
0             & 0             & 0             & 0             & 0             &  2\sqrt{14}   & 0             & -2\sqrt{2}    & 0              \\
0             & -6            & 0             &  2\sqrt{7}    & 0             & 0             & 0             & 0             & 0              \\
0             & 0             & 0             & 0             & 0             &  2\sqrt{2}    & 0             &  2\sqrt{14}   & 0              \\
0             &  2\sqrt{7}    & 0             &  6            & 0             & 0             & 0             & 0             & 0              \\
0             & 0             & 0             & 0             &  3            & 0             &  2\sqrt{5}    & 0             &  \sqrt{35}     \\
-2\sqrt{14}   & 0             & -2\sqrt{2}    & 0             & 0             & 0             & 0             & 0             & 0              \\
0             & 0             & 0             & 0             &  2\sqrt{5}    & 0             &  4            & 0             & -2\sqrt{7}     \\
 2\sqrt{2}    & 0             & -2\sqrt{14}   & 0             & 0             & 0             & 0             & 0             & 0              \\
0             & 0             & 0             & 0             &  \sqrt{35}    & 0             & -2\sqrt{7}    & 0             &  1
\end{bmatrix},
\end{equation*}
\end{Small}

\begin{equation*}
R_y(\beta) = R_x(\frac{\pi}{2}) \times R_z(\beta) \times R_x(\frac{\pi}{2})^T,
\end{equation*}

\begin{equation*}
R_x(\alpha) = R_y(\frac{\pi}{2})^T \times R_z(\alpha) \times R_y(\frac{\pi}{2}).
\end{equation*}

See \cite{Blanco1997,Choi1999,Collado1989,Ivanic1996} for more details.

\newpage

\section{Appendix A.2}
The rotational matrices $R_x$, $R_y$ and $R_z$ for spherical harmonics of degree 3 (octupoles) are defined as follows:

\begin{Small}
\begin{equation*}
R_z(\gamma) =
\begin{bmatrix}
\cos 3\gamma  & 0             & 0             & 0             & 0             & 0             &  \sin 3\gamma \\
0             & \cos 2\gamma  & 0             & 0             & 0             &  \sin 2\gamma & 0             \\
0             & 0             & \cos  \gamma  & 0             &  \sin  \gamma & 0             & 0             \\
0             & 0             & 0             & 1             & 0             & 0             & 0             \\
0             & 0             & -\sin  \gamma & 0             & \cos  \gamma  & 0             & 0             \\
0             & -\sin 2\gamma & 0             & 0             & 0             & \cos 2\gamma  & 0             \\
-\sin 3\gamma & 0             & 0             & 0             & 0             & 0             & \cos 3\gamma                
\end{bmatrix},
\end{equation*}

\begin{equation*}
R_x(\frac{\pi}{2}) = \frac{1}{4}
\begin{bmatrix}
0             & 0             & 0             & \sqrt{10}     & 0             & -\sqrt{6}     & 0             \\
0             & -4            & 0             & 0             & 0             & 0             & 0             \\
0             & 0             & 0             & \sqrt{6}      & 0             & \sqrt{10}     & 0             \\
-\sqrt{10}    & 0             & -\sqrt{6}     & 0             & 0             & 0             & 0             \\
0             & 0             & 0             & 0             & -1            & 0             & -\sqrt{15}    \\
\sqrt{6}      & 0             & -\sqrt{10}    & 0             & 0             & 0             & 0             \\
0             & 0             & 0             & 0             & -\sqrt{15}    & 0             & 1             
\end{bmatrix},
\end{equation*}
\end{Small}

\begin{equation*}
R_y(\beta) = R_x(\frac{\pi}{2}) \times R_z(\beta) \times R_x(\frac{\pi}{2})^T,
\end{equation*}

\begin{equation*}
R_x(\alpha) = R_y(\frac{\pi}{2})^T \times R_z(\alpha) \times R_y(\frac{\pi}{2}).
\end{equation*}

See \cite{Blanco1997,Choi1999,Collado1989,Ivanic1996} for more details.

\bibliographystyle{plain}
\bibliography{lit}

\end{document}